\documentclass[a4paper]{amsart}

\usepackage{amsfonts}
\usepackage{amssymb}
\usepackage{amsthm}
\usepackage{amsmath}
\usepackage{pst-node}

\usepackage{tikz-cd}
\usepackage{pstricks}
\usepackage{newlfont}
\usepackage{graphicx}
\usepackage[all,arc]{xy}
\usepackage{epsfig}
\usepackage{amsmath,amscd}
\usepackage{url}
\usepackage{todonotes}
\usepackage{leftidx}
\usepackage{multicol}
\usepackage[T1]{fontenc}
\usepackage{mathtools}
\usepackage[noend]{algpseudocode}
\usepackage[section]{algorithm}
\usepackage{subfigure}
\newcommand{\dst}{\displaystyle}

\definecolor{ududff}{rgb}{0.30196078431372547,0.30196078431372547,1.}
\definecolor{cqcqcq}{rgb}{0.7529411764705882,0.7529411764705882,0.7529411764705882}

\usepackage{pgf}
\usepackage{tikz}
\usepackage[pagewise]{lineno}
\usetikzlibrary{arrows,automata}

\newcommand {\dz}  {\mbox{\bf Proof: }}

\newtheorem{proposition}{Proposition}[section]
\newtheorem {theorem}{Theorem}[section]

\newtheorem {definition}{Definition}[section]

\newtheorem{remark}{Remark}

\usepackage[utf8]{inputenc}

\title{Homology of polyomino tilings on flat surfaces}
\author{Edin Li\dj{}an and Djordje Barali\'{c}}

\address{\scriptsize{Faculty of Pedagogy, University of Biha\'{c}, Biha\'{c},  Bosnia and Herzegovina}}
\email{lidjan\_edin@hotmail.com}
\address{ \scriptsize{Mathematical Institute SASA, Belgrade, Serbia }}
\email{djbaralic@mi.sanu.ac.rs}


\begin{document}

\maketitle

\begin{abstract}
The homology group of a tiling introduced by M. Reid is studied for certain topological tilings. As in the planar case, for finite square grids on topological surfaces, the method of homology groups, namely the non-triviality of some specific element in the group allows a `coloring proof' of impossibility of a tiling. Several results about the non-existence of
polyomino tilings on certain square-tiled surfaces are proved in the paper.
\end{abstract}

\section{Introduction}

Recreational mathematics comprises various subjects including combinatorial games, puzzles, card tricks, art, etc.  Its problems are typically easily understood by a general audience, yet their solution
often requires rigorous research. Indeed, a significant number of mathematical disciplines have been
grounded on ideas sparked by challenges from recreational mathematics. For example, graph theory
has its roots in the solution of the problem of The Seven Bridges of K\"{o}nigsberg, and magic squares
contributed to the foundations of combinatorial designs.

\begin{figure}[h]
  \centering
    \includegraphics[width=0.5\textwidth]{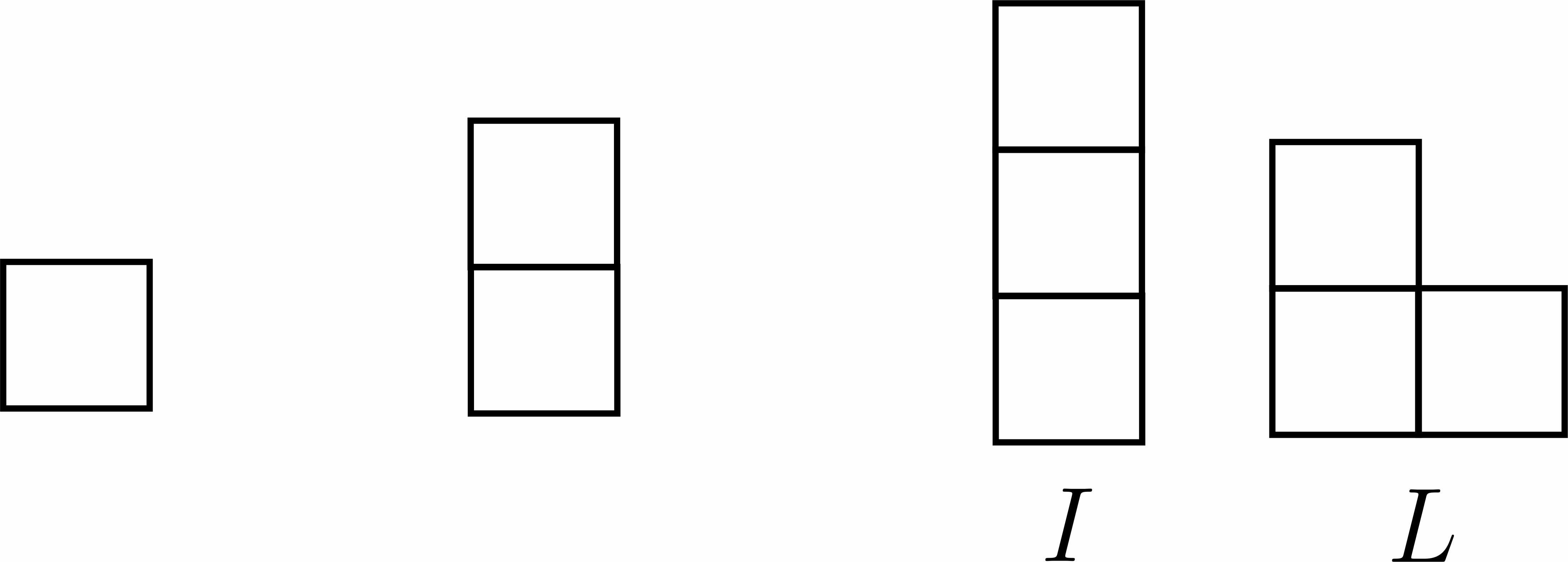}
    \caption{Monomino, domino and trominoes}
    \label{domino}
  \end{figure}

\begin{figure}[h]
    \centering
        \includegraphics[width=0.7\textwidth]{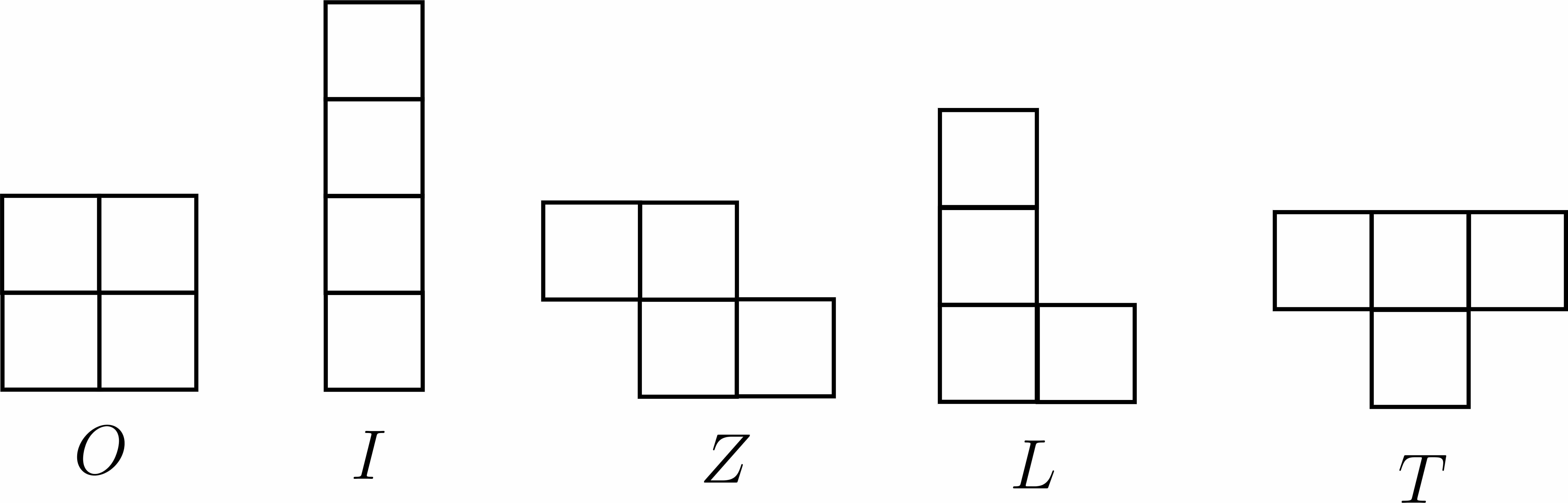}
    \caption{Tetrominoes}
    \label{tetramino}
    \end{figure}

\begin{figure}[!ht]
  \centering
    \includegraphics[width=\textwidth]{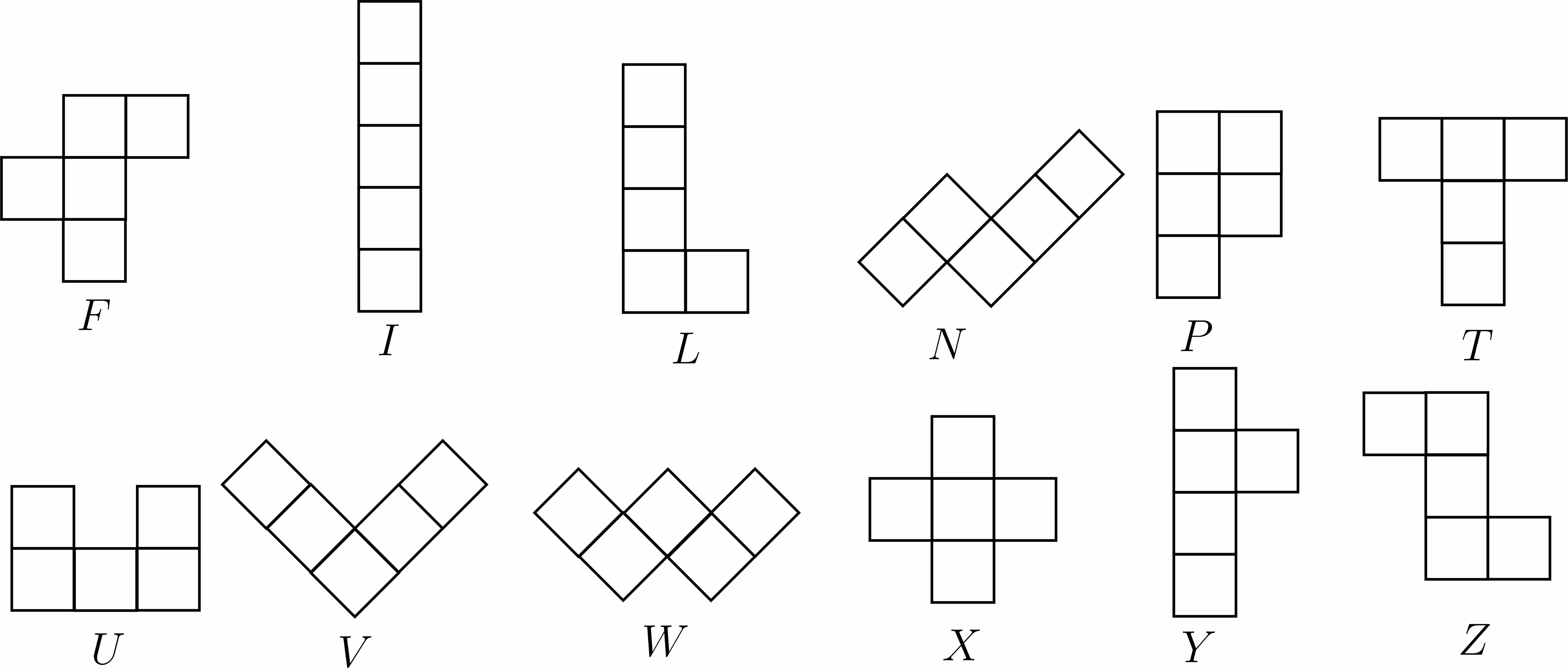}
    \caption{Pentominoes}
    \label{pentamino}
\end{figure}

A \textit{polyomino} is a planar geometric figure formed by
joining one or more identical squares edge-to-edge. It may also be regarded as a finite subset of the regular square grid with a connected interior. A polyomino consisting of exactly $n$ cells is
called an \textit{$n$-omino}.  Polyomino shapes for $n\leq 5$ are
illustrated in Figures \ref{domino}, \ref{tetramino} and
\ref{pentamino}. Some polyominoes were named after letters of the alphabet closely resembling
them, as can be seen in Figures \ref{tetramino}  and \ref{pentamino}. They were popularized by Salmon Golomb who wrote the first monograph on
polyominoes \cite{golomb1}, and by Martin Gardner in his Scientific
American columns ``Mathematical Games'', see \cite{gardner}. In fact, the word polyomino was coined by Golomb in \cite{GolombM}. Today they
are one of the most popular subjects of recreational mathematics, being of great interest to not only
mathematicians but physicists, biologists, and computer scientists as well. For more information,
we refer the reader to surveys \cite{ardila} and
\cite{handbook}.

\textit{The polyomino tiling problem} asks whether it is possible to properly tessellate a finite region of cells,
say $M$, with polyomino shapes from a given set $\mathcal{T}$. There are numerous generalizations of this question
for symmetric and asymmetric tilings, higher dimensional analogs, polyomino type problems on other regular lattice grids (triangular, hexagonal), etc. However, the problem is NP-hard in general, and we
can give definite answers only in a limited number of cases.

This enthralling problem from recreational mathematics has attracted attention of both mathematicians and non-experts. There were many results establishing criteria for proper tilings by some specific polyomino shapes (see \cite{golomb3}, \cite{golomb4}, \cite{golomb5}, \cite{reid2}
and \cite{reid3}).  Conway and Lagarias developed in \cite{cl} the
so-called `boundary-word method' for addressing this
question. Their ideas were further developed by Reid in \cite{reid}
who  assigned to each set of tiles $\mathcal{T}$ \textit{the
homology and the homotopy group of tilings} and formulated a
necessary condition for existence of a proper tiling of a finite
region $M$ in a plane.

Reid's powerful idea allows natural generalization to a much bigger class of combinatorial tilings. Instead of considering planar regions, we study regions which are obtained by identifying parts of the boundary of a planar region resulting in a flat Riemann surface. The only flat compact Riemann surfaces are the torus and Klein bottle, but one can give higher-genus surfaces a flat metric everywhere except at certain cone points, and then remove neighborhoods of the singular points to get a flat surface with boundary. Surfaces with a flat metric obtained by pairwise identification of sides of a collection of plane polygons via translations of their sides, are called translation surfaces. Translation surfaces can also be defined as Riemann surfaces with a holomorphic 1-form. In particular, we are interested in a
subclass of translation surfaces called a square-tiled surface.  A square-tiled surface is any translation
surface obtained from a polygon $P$ which is itself obtained by putting a collection of copies of the unit square side by side. In general, the total angle around a corner of a square of a square-tiled surface $S$ is
a non-trivial multiple of $2\pi$. Any such point is called a conical singularity of $S$. In this paper, we study
the problem of tiling a surface S subdivided into a finite `combinatorial' grid by a finite set of
polyomino shapes T and define the homology group $H_S (\mathcal{T})$.

Square-tiled and translation surfaces arise in dynamical systems, where they can be used to model billiards, and in Teichm\"{u}ller theory. They have a rich mathematical structure and may be studied from multiple points of view (flat geometry, algebraic geometry, combinatorial group theory, etc.). We present some new results and illustrate examples explaining the application of the homology group of generalized polyomino type tilings in the combinatorial and the topological context.

In Section 2, we introduce the homology tiling group for finite square grids on surfaces with boundaries based on \cite{reid}.
Several results about the impossibility of tiling certain concrete square-tiled surfaces are proved using the homology group of the tiling in Section 3. Our main novelty lies in Theorem \ref{teo3.7} which establishes a general result connecting the I-polyomino shape with the genus of the surface.

\section{Tiling Problem on Surfaces}

\begin{figure}[h]
  \centering
    \includegraphics[width=0.7\textwidth]{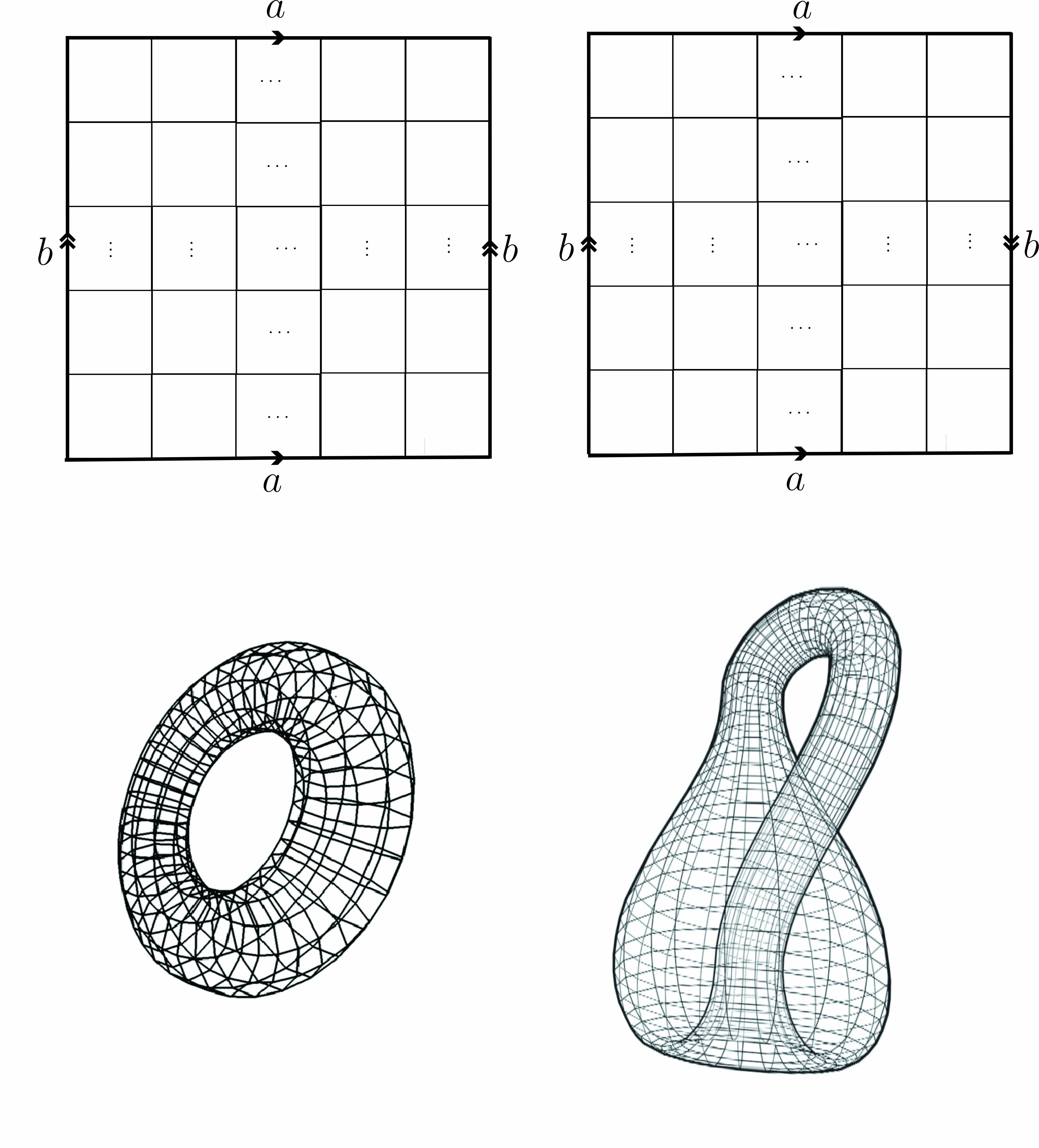}
    \caption{Square grids on a torus and a Klein bottle}
    \label{torus}
\end{figure}

The standard square grid in the plane is characterized by the property that exactly four edges meet at each vertex, and each vertex is shared by four squares in the grid. We assume that every edge in a combinatorial grid on a surface is shared by exactly two squares, unless it is on the boundary. This local property allows us to define a polyomino tiling on a topological surface in the same way as in the planar case, and we will refer to such a structure as \textit{the square grid on a surface}. For example, identification of parallel edges of the boundary of $m\times n$ grid in the same directions provides such a grid on torus. Identification of two pairs of parallel sides of an $m\times m$, $m\geq 3$ square, but in the opposite direction in one of the pairs,  provides examples of square grids on Klein bottle, see Figure \ref{torus}. However, if the surface has no boundary, then each
vertex is shared by four squares, so the number of vertices equals the number of squares. Likewise,
each edge is shared by two squares, so the number of edges is twice the number of squares. This makes
the Euler characteristic
\begin{eqnarray*}
\chi (M) = V- E + F = F - 2F + F = 0,
\end{eqnarray*}
so $M$ is either the torus or the Klein bottle.

On topological surfaces with boundaries, square grids are not rare structures. One way to obtain them is
by identification of certain faces of a finite region in a planar square grid. Identification of faces allows
the additional possibility for placing a polyomino tile, so we have to develop means to treat tiling
problems. Surfaces obtained by gluing sides of a polygon are extensively studied in mathematics and
this is an interesting research topic in itself (see \cite{Ahmedov}, \cite{harer},
\cite{Koch} and \cite{stillwel}).

Actually, above mentioned combinatorial structures are directly related to mathematical concepts known as \textit{translation surfaces}. Combinatorially, a translation surface may be defined in the following way. Let $\displaystyle P_{1},\ldots ,P_{m}$ be a collection of  polygons in the Euclidean plane and suppose that for every side $\displaystyle s_{i} $ of any $\displaystyle P_{k}$ there is a side $\displaystyle s_{j}$ of some $\displaystyle P_{l}$ with $\displaystyle j\not =i$ and $\displaystyle s_{j}=s_{i}+{\vec {v}}_{i}$ for some nonzero vector $\displaystyle {\vec {v}}_{i}$ and so that $\displaystyle {\vec {v}}_{j}=-{\vec {v}}_{i}$. The space obtained by identifying all $\displaystyle s_{i}$ with their corresponding $\displaystyle s_{j}$ through the map $\displaystyle x\mapsto x+{\vec {v}}_{i}$ is a translation surface.

A particular class of translation surfaces known as square-tiled surfaces is of wide interest for mathematics. A square-tiled surface is an orientable connected surface obtained from a finite collection of unit squares in a plane after identifications of pairs of parallel sides via adequate translations. In general, the total angle around a corner of a square of a square-tiled surface $M$ is a non-trivial multiple of $2\pi$ and any such point is called a conical singularity of $M$. In our considerations we will consider flat surfaces with cone points with cone angle a multiple of $\frac{\pi}{2}$.

The tiling problem for a finite subset of the regular planar square grid by a finite set of polyomino
prototiles has been studied extensively in the past few decades. However, there exist many other
topological $2$-manifolds which admit subdivision into a finite number of squares which preserves the
structure of regular square grid and for which the tiling problem is also defined. One natural way to
obtain such structures is by gluing some of the faces of a finite subset of regular square lattice in the
plane, and some results and examples of polyomino tiling problems in this context are known in
literature under the notion of topological tilings. Special cases of cylinders, torus, M\"{o}bius strip, Klein bottle and projective plane with a $2$-disk removed were studied in \cite{golomb1}, \cite{remila} and \cite{lima}.

Several techniques for finding obstructions to tiling are known, and one of the most charming is that of
a `generalized chessboard coloring'. This method rests on the fact that the chessboard with two
opposite square corners removed cannot be tiled by dominoes, as the difference between the number of
white and black squares is two, see \cite{GolombM}. The general idea is to use several colors and color the squares
of the considered region in a special pattern `sensitive' to the given set of polyominoes. In other words,
the coloring imposes some number theoretical condition which serves as an obstruction to a tiling.
However, it is not easy to find a coloring argument for proving nonexistence of a tiling. Michael Reid
introduced in \cite{reid} the so-called homology group of a tiling and showed that proof of nontriviality of a
special element in this group assigned to the finite subset of regular square lattice produces a
generalized chessboard coloring argument. His homology tiling group method is therefore at least as
powerful as the coloring argument. In the same paper, Reid gave many examples where the tiling
homology group is inefficient for proving non-existence of a tiling.

The problem of polyomino tilings was studied by Conway and
Lagarias in \cite{cl} where they introduced a
new technique using boundary word invariants to formulate necessary conditions for the existence of
tilings. Based on their ideas, Reid presented in
\cite{reid} a new strategy for treating tiling problems, working
with the so-called homotopy group of tiling. Reid's \textit{homotopy tiling group method} was so far the most
successful in establishing necessary criteria for existence of tilings.

Our main observation is that Reid's tiling homology group method can be applied to studying
topological tilings. A standard model for obtaining topological surfaces is identification of sides of a
polygon and as clearly presented in \cite{Koch} and
\cite{stillwel}.

Let $M$ be a topological surface with boundary obtained by gluing
of sides of some finite subset $R$ of the regular square grid in the
plane and let $\mathcal{T}$ be a finite set of polyomino tiles.
Gluing of faces provides more ways for placement of tiles from
$\mathcal{T}$ on $M$ then in the case of $R$, so $M$ may be tiled
even if $R$ does not admit a tiling by tiles from $\mathcal{T}$.
We introduce the tiling homology group $H(M, \mathcal{T})$ in the same
fashion  as Michael Reid.

Let $A$ be the free abelian group generated by the set of  cells
of $M$. We assume that all cells of $M$ preserve labeling by $(i,
j)$ from $R$. The generator of $A$ corresponding to the cell $(i,
j)$ is denoted by $a_{i, j}$. Let $B(M, \mathcal{T})$ be the
subgroup generated by elements  corresponding to all possible
placements of tiles in $\mathcal{T}$, i.e. by the sums of elements
assigned to cells of $M$ that can be covered by a tile from
$\mathcal{T}$.

\begin{definition}
The tiling homology group of $(M, \mathcal{T})$ is the quotient
group  $$H(M, \mathcal{T}) = A/B(M, \mathcal{T}).$$
\end{definition}

Let us denote by $\bar{a}_{i, j}$ the image of $a_{i, j}$ in $H(M,
\mathcal{T})$. As in the planar case, there is an element $\Theta \in
H(M, \mathcal{T})$ assigned to $M$
$$\Theta:= \sum_{(i, j)\in M} \bar{a}_{i, j}$$ which is clearly zero when there is a tiling of $M$ by polyominoes from $\mathcal{T}$. Thus, $\Theta$ is an obstruction to tiling. Recall that Reid considered in his paper the so-called \textit{signed tiling}, where he allowed polyomino tiles to have positive and negative signs. Clearly, the signed tiling of $M$ by $\mathcal{T}$ exists if and only if $\Theta$ is trivial in $H(M, \mathcal{T})$.

Reid's \cite[Proposition~2.10]{reid} also holds for topological
tilings by polyominoes. It states that nontrivial $\Theta$
produces special numbering of cells in $M$ that yields a
generalized chessboard coloring argument. We adapt his proof to
the case of topological tilings.

\begin{proposition}
Let $M$ be a topological surface with boundary with a finite
square grid and finite set of polyominoes $\mathcal{T}$ such that
$\Theta$ is nontrivial in $H(M, \mathcal{T})$. Then there is the
numbering of the cells in $M$ by rational numbers such that
\begin{itemize}
\item[i)] for any placement of a tile from $\mathcal{T}$, the
total sum of covered numbers is an integer, and \item[ii)]  the
total covered by the cells of $M$ is not an integer.
\end{itemize}
\end{proposition}
\noindent \dz Consider the cyclic subgroup $\langle \Theta \rangle
\subset H(\mathcal{T})$ generated by $\Theta$. We define a
homomorphism $\varphi: \langle \Theta \rangle \rightarrow
\mathbb{Q}/\mathbb{Z}$ with $\varphi(\Theta)\neq 0$. If $\Theta$
has infinite order we set $\varphi(\Theta) = \frac{1}{2}
\mod{\mathbb{Z}}$, while if $\Theta$ has finite order $n > 1$,
then we define $\varphi(\Theta) = \frac{1}{n} \mod {\mathbb{Z}}$.
Since $\mathbb{Q}/\mathbb{Z}$ is a divisible abelian group, the
homomorphism $\varphi$ extends to a homomorphism $H(M,
\mathcal{T}) \rightarrow \mathbb{Q}/\mathbb{Z}$, also called
$\varphi$. Here we used the familiar fact about equivalence of the notions of injective group and divisible group for abelian groups \cite[Proposition~6.2]{bredon}. Since $A$ is a free abelian group, the composite map
\[\begin{tikzcd}A \arrow[two heads]{r} & A/B(M, \mathcal{T}) =
H(M, \mathcal{T}) \arrow{r}{\varphi} &
\mathbb{Q}/\mathbb{Z}\end{tikzcd}
\] lifts to a homomorphism $\psi: A \rightarrow\mathbb{Q}$, such that the following diagram commutes

\[\begin{tikzcd}
A \arrow{r}{\psi} \arrow[two heads]{d} & \mathbb{Q} \arrow[two heads]{d} \\
H(M, \mathcal{T}) \arrow{r}{\varphi} & \mathbb{Q}/\mathbb{Z}
\end{tikzcd}
\]
where the vertical surjections are the quotient maps. Desired
numbering of the cells is defined by $\psi$, and since $B(M,
\mathcal{T})$ is in the kernel of $A \rightarrow
\mathbb{Q}/\mathbb{Z}$, every tile placement covers an integral
total. But, $\varphi({\Theta})\neq 0$ and total of the cells in
$M$ is not an integer. \qed

Reid's tiling homology group was systematically studied using Gr\"{o}bner bases in the works of Muzika-Dizdarevi\'{c}, Timotijevi\'{c} and \v{Z}ivaljevi\'{c}, see \cite{muzziv} and \cite{muz}.

\section{Nonexistence of polyomino tilings on surfaces}

In this section we prove several results on nonexistence of tilings
on surfaces of different genus with boundaries by some given
polyomino sets as an illustration of the homology method.

First we formulate three results for polyomino tilings on a torus
square grid. Such cases were also studied in the past
\cite{remila} as they are close to the planar case.

\begin{figure}[H]
    \centering
    \includegraphics[width=0.60\linewidth]{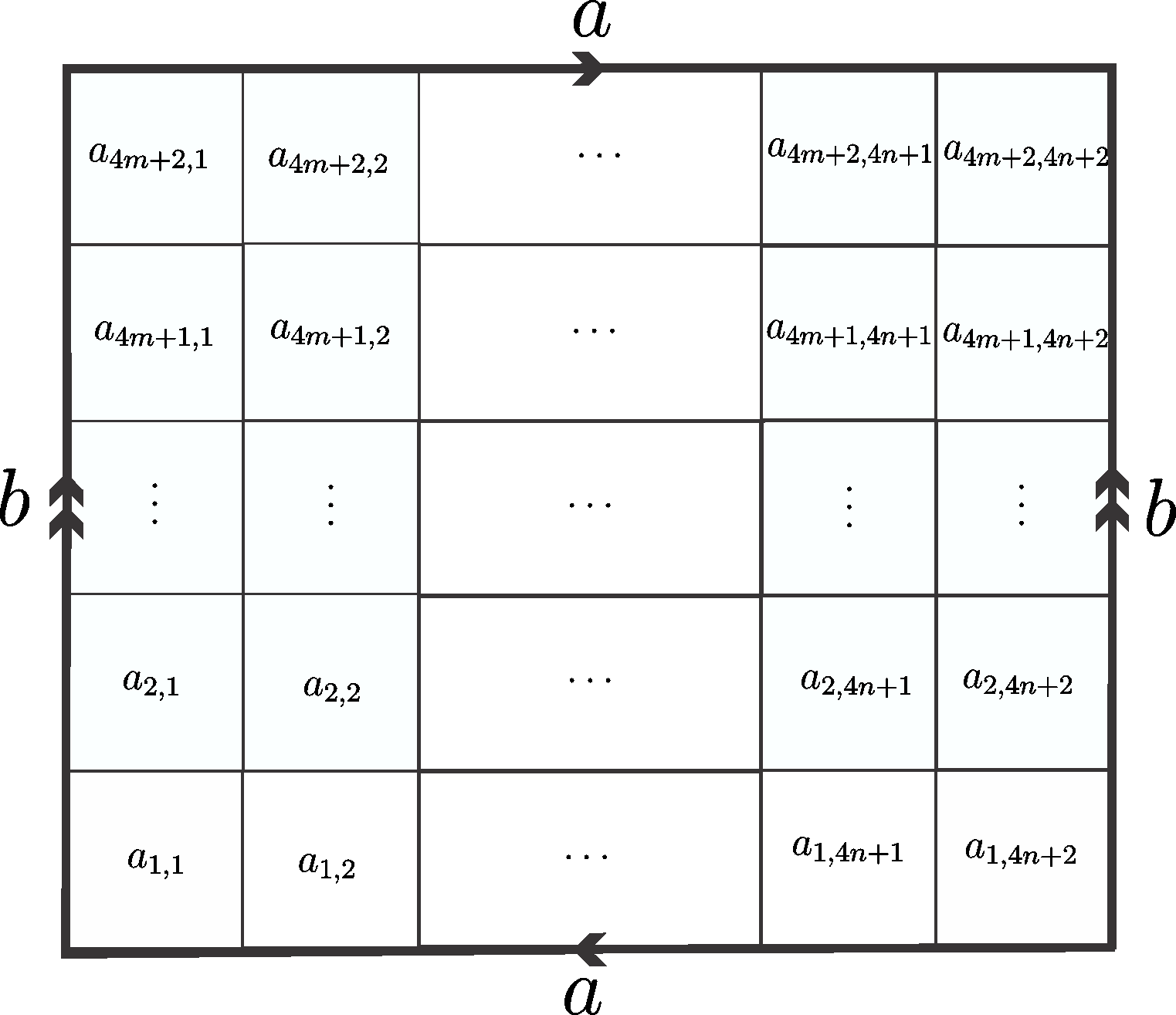}
    \caption{Torus grid of dimension $(4m+2) \times (4n+2)$}\label{I_mreza}
\end{figure}

\begin{theorem}\label{teo1}
A square torus grid of dimension  $(4m+2)\times (4n+2)$ cannot be
tiled by I-tetrominoes, see Figure \ref{tetramino}.
\end{theorem}
\noindent \dz Consider a  $(4m+2) \times (4n+2)$ square torus grid
model in a plane with cells labelled as in Figure  \ref{I_mreza}.

Investigate all possible placements of a tile in the given model.
To each placement one can assign one of two types of  relations:
\begin{eqnarray*}
\bar{a}_{i, j} + \bar{a}_{i,j+1} + \bar{a}_{i,j+2} + \bar{a}_{i,
j+3} =0 \quad \text{and}\quad \bar{a}_{i, j} + \bar{a}_{i+1,j} +
\bar{a}_{i+2,j} + \bar{a}_{i+3, j} =0
\end{eqnarray*}
where $i=1,\ldots, 4m+2$ labels a row, and $j=1, \ldots, 4n+2$
labels a column on the given torus grid. We assume that indices of
rows is modulo $4m+2$ and modulo $4n+2$ for columns in the
relations above. Considering the relation
\begin{eqnarray*} \bar{a}_{i, j+1} + \bar{a}_{i,j+2} + \bar{a}_{i,j+3} + \bar{a}_{i, j+4}
=0
\end{eqnarray*}  we obtain that in  the homology group of this tiling it holds that
\begin{eqnarray*}
\bar{a}_{i,j} = \bar{a}_{i,j+4}
\end{eqnarray*} for all $i, j\in \{1, 2,\dots, 4k+2\}$.
Analogously,  $\bar{a}_{i, j}=\bar{a}_{i+4,j}$.

\begin{figure}[H]
    \centering
    \includegraphics[width=0.60\linewidth]{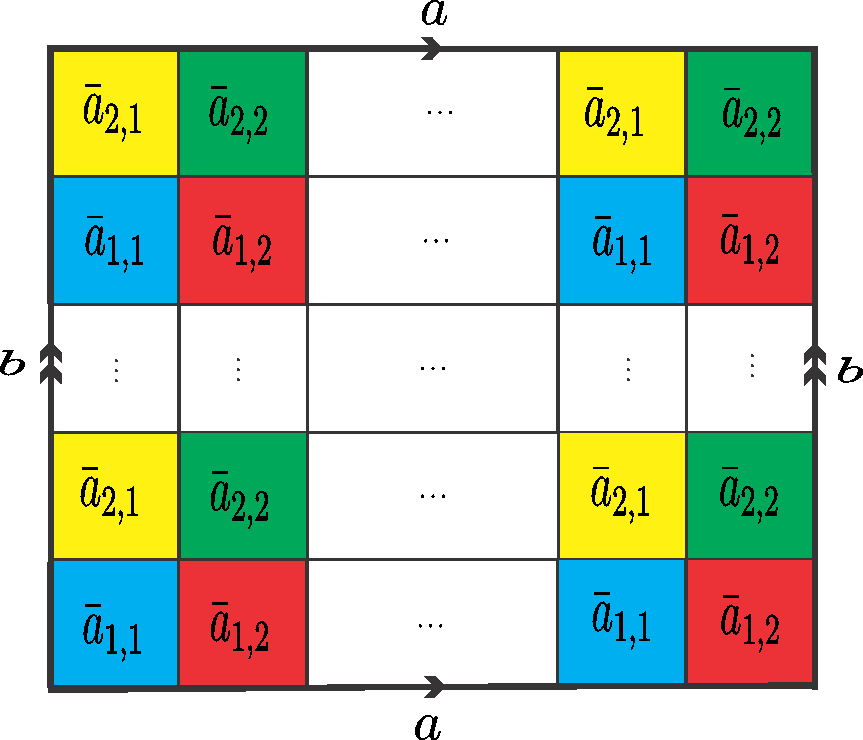}
    \caption{Coloring of the equivalent cells of the torus grid}\label{T36_4}
\end{figure}

From the relations corresponding to placements over the identified
faces of the rectangle representing our torus grid, we obtain
additional cells of the grid whose corresponding generators in the
homology group of tiling are equal.  Using
\begin{eqnarray*}
\bar{a}_{i, 4m-1} + \bar{a}_{i,4m} + \bar{a}_{i,4m+1} + \bar{a}_{i, 4m+2} =0 \quad \text{and} \\
\bar{a}_{i, 4m} + \bar{a}_{i,4m+1} + \bar{a}_{i,4m+2} +
\bar{a}_{i, 1} =0
\end{eqnarray*}
we conclude that $\bar{a}_{i, 1}=\bar{a}_{i,4m-1}$. In the same
fashion we deduce that $\bar{a}_{i, 2}=\bar{a}_{i,4m}$,
$\bar{a}_{1, i}=\bar{a}_{4n-1,i}$ and $\bar{a}_{2,
i}=\bar{a}_{4n,i}$ for all $i$. Combining the equalities above, we
obtain
\begin{eqnarray*}
\bar{a}_{i, j}=\begin{cases}
\bar{a}_{1, 1}, \quad \text{if} \quad i \equiv 1 \pmod{2}, \, j \equiv 1 \pmod{2},\\
\bar{a}_{1, 2}, \quad \text{if} \quad i \equiv 1 \pmod{2}, \, j \equiv 0 \pmod{2},\\
\bar{a}_{2, 1}, \quad \text{if} \quad i \equiv 0 \pmod{2}, \, j \equiv 1 \pmod{2},\\
\bar{a}_{2, 2}, \quad \text{if}  \quad i \equiv 0 \pmod{2}, \, j \equiv 0 \pmod{2}.\\
\end{cases}
\end{eqnarray*} as depicted in Figure \ref{T36_4}.

If we put I-tetromino shape on the torus grid with equivalent
cells we obtain one of the following relations
\begin{align*}
        2\bar{a}_{1,1} + 2\bar{a}_{1,2} & =  0, & 2\bar{a}_{1,1} + 2\bar{a}_{2,1} & =  0,  & \\
        2\bar{a}_{2,1} + 2\bar{a}_{2,2} & =  0, & 2\bar{a}_{1,2} + 2\bar{a}_{2,2} & =  0.   &
\end{align*}
Therefore,  our homology group is isomorphic to the quotient group
of the free abelian group with four generators  by the four
relations given above. Let us observe that one of these relations can be obtained from the remaining three so we can omit the relation $2\bar{a}_{2,1} + 2\bar{a}_{2,2} =  0$ . We can consider the presentation of the group  using the
following four generators $a=\bar{a}_{1,1}$, $b=\bar{a}_{1,1} +
\bar{a}_{1,2}$, $c=\bar{a}_{1,1} + \bar{a}_{2,1}$ and
$d=\bar{a}_{2,2} -\bar{a}_{1,1}$. It is clear that  $\dst
2b=2c=0$, and little more effort gives $\dst 2d =0$. Thus, our homology group of tiling is isomorphic to
 $$\dst G \langle a, b, c, d |2b=2c=2d=0 \rangle\cong \mathbb{Z}\oplus (\mathbb{Z}_2)^3. $$

It is easily seen that everything but
the top two (or bottom two) rows of our grid are easily tiled by vertical
I-tetrominoes, and that in the top two rows everything but the right-most two
columns are tiled by horizontal I-tetrominoes, so $\Theta$ is the sum of elements corresponding to the four upper right cells. Thus,
\begin{align*}\label{jed1}
\Theta & =
\bar{a}_{1,1}+\bar{a}_{1,2}+\bar{a}_{2,1} + \bar{a}_{2,2}=b+c+d
\end{align*}
is nontrivial in the tiling homology group, so desired tiling is not
possible. \qed

\begin{remark} We can reach the same conclusion using coloring of the square torus grid as in Figure \ref{T36_4}. Each tile covers  $2$ blue and  $2$ yellow cells, or $2$ blue and $2$ red, or $2$ yellow and $2$ green, or $2$ red and $2$ green. Since the number of cells of  each color is odd and each  tile covers an even number of cells of the same color, we conclude that tiling is not possible.
\end{remark}

\begin{figure}[H]
    \centering
    \includegraphics[width=0.60\linewidth]{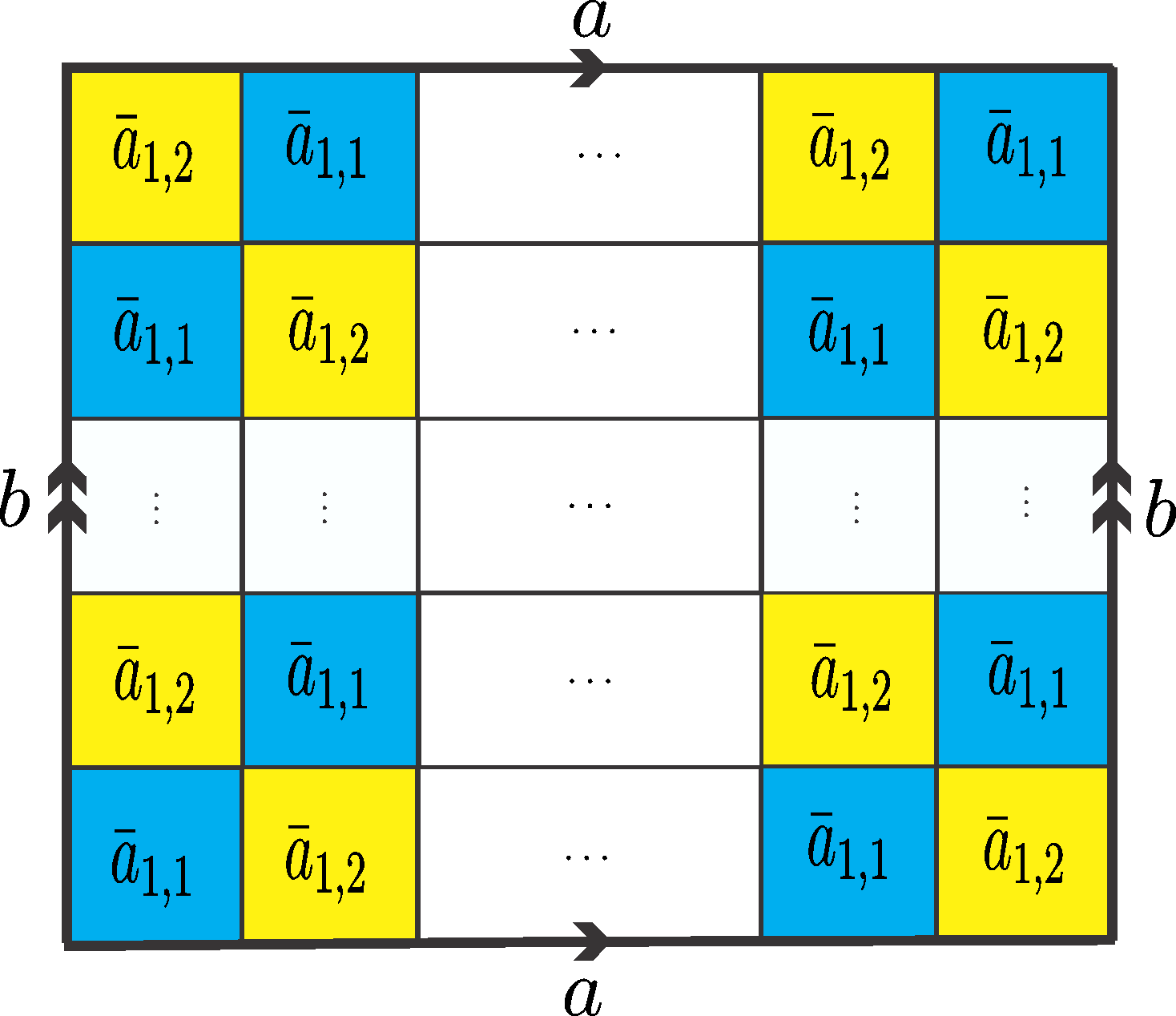}
    \caption{Coloring of the equivalent cells of the torus grid}\label{Teorema_2_3}
\end{figure}

\begin{theorem}
A square torus grid of dimension $(4m+2)\times (4n+2)$ cannot be
tiled with $T$ tetrominoes.
\end{theorem}
\noindent \dz Consider the torus grid presented as in  Figure
\ref{I_mreza}. Consider all possible placements of T tetromino. To
each placement we can assign one of the following relations:
\begin{eqnarray*}\label{Tjed1}
\bar{a}_{i, j} + \bar{a}_{i,j+1} + \bar{a}_{i,j+2} + \bar{a}_{i+1, j+1} &=&0,  \\
\bar{a}_{i, j} + \bar{a}_{i,j+1} + \bar{a}_{i,j+2} + \bar{a}_{i-1, j+1} &=&0,
 \\
\bar{a}_{i,j}+\bar{a}_{i+1,j}+\bar{a}_{i+2,j}+\bar{a}_{i+1,j+
1}&=&0 \quad \mbox{and} \\
\bar{a}_{i,j}+\bar{a}_{i+1,j}+\bar{a}_{i+2,j}+\bar{a}_{i+1,j-
1}&=&0,
\end{eqnarray*}
where we use the same labelling as in the proof of Theorem
\ref{teo1}. From them we directly deduce that  in the  homology
group of tiling it holds that
\begin{eqnarray*}
    \bar{a}_{i+2,j}=\bar{a}_{i,j} = \bar{a}_{i,j+2}
\end{eqnarray*}
for all $i$ and $j$. Therefore,  \begin{eqnarray*} \bar{a}_{i,
j}=\begin{cases}
\bar{a}_{1, 1}, \quad \text{if} \quad i-j\equiv 0 \pmod{2}, \\
\bar{a}_{1, 2}, \quad \text{if} \quad i-j \equiv 1 \pmod{2},
\end{cases}
\end{eqnarray*}
as it is illustrated in Figure \ref{Teorema_2_3}.

Placing a T-tetromino shape on the torus grid with equivalent
cells, we get one of the following two relations
\begin{eqnarray*}
        3\bar{a}_{1,1} + \bar{a}_{1,2} & = & 0 \quad \mbox{and} \\
        3\bar{a}_{1,2} + \bar{a}_{1,1} & = & 0.
\end{eqnarray*}
Therefore, our homology group is isomorphic to the group
$$\dst G \langle \bar{a}_{1, 1} |8\bar{a}_{1,1}  = 0 \rangle \cong \mathbb{Z}_{8}.$$ Our grid has $2m$ cells $\bar{a}_{1,1}$ and $\bar{a}_{1,2}$, where $k=(2m+1)(2n+1)$. So the element that corresponds to this grid
\begin{eqnarray*}
\Theta &  = & 2k \bar{a}_{1,1} + 2k \bar{a}_{1,2}  =  -4k
\bar{a}_{1,1}= 4 \bar{a}_{1,1}
\end{eqnarray*}
is nontrivial in the homology group, so desired tiling does not
exist. \qed

\begin{remark}
The same conclusion can be obtained using  coloring in Figure
\ref{Teorema_2_3} and  parity argument for the total number of
cells in the grid.
\end{remark}

\begin{theorem}
A square torus grid of dimension $(4m+2)\times (4n+2)$ cannot be
tiled with $X$ hexominoes (Figure \ref{produzeni_krst}).

\begin{figure}[H]
    \centering
    \includegraphics[width=0.20\linewidth]{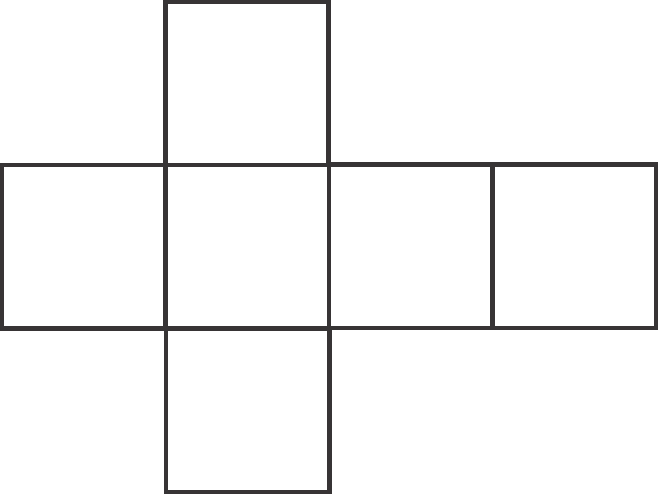}
    \caption{$X$ hexomino}\label{produzeni_krst}
\end{figure}

\end{theorem}
\noindent \dz Consider planar model of torus grid of dimension
$(4m+2) \times (4n+2)$  as in Figure  \ref{I_mreza}. Examine all
possible horizontal placements of our tile. Each of them yields a
relation
\begin{eqnarray}\label{jednacinareda}
\bar{a}_{i, j} + \bar{a}_{i,j+1} + \bar{a}_{i,j+2} + \bar{a}_{i,
j+3} + \bar{a}_{i+1,j+1} + \bar{a}_{i-1,j+1} =0,
\end{eqnarray}
where where the rows and columns are labelled analogously as in
the proof of Theorem \ref{teo1}. From (\ref{jednacinareda}) we
conclude that in the homology group of this tilling  it holds
$\bar{a}_{i,j}=\bar{a}_{i,j+4}$ for all $i$ and $j$. Since
$\bar{a}_{i,4n-1}=\bar{a}_{i,1}$, $\bar{a}_{i,4n}=\bar{a}_{i,2}$,
$\bar{a}_{i,4n+1}=\bar{a}_{i,3}$ and
$\bar{a}_{i,4n+2}=\bar{a}_{i,4}$ we further get that for all  $i$
and $j$ it also holds $\bar{a}_{i,j}=\bar{a}_{i,j+2}$.

\begin{figure}[h!]
    \centering
    \includegraphics[width=0.60\linewidth]{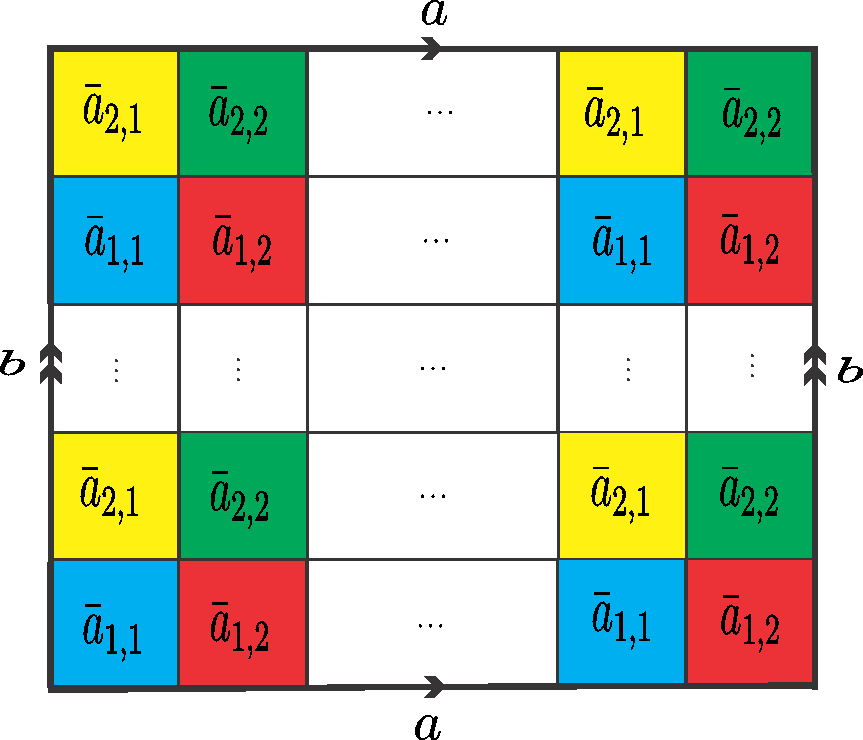}
    \caption{Coloring of the equivalent cells of the torus grid}\label{Teorema_3_3}
\end{figure}

Analogous consideration of vertical placements implies
$\bar{a}_{i,j}=\bar{a}_{i+2,j}$ for all $i$ and $j$. Equivalences
of the cells in the grid in the homology group of tiling are
depicted in Figure \ref{Teorema_3_3}.

Thus, we deduce that the homology group of tiling is the quotient
of the free abelian group with four generators $\dst G \langle
\bar{a}_{1,1}, \bar{a}_{1,2}, \bar{a}_{2,1}, \bar{a}_{2,2}\rangle
$ modulo following relations
\begin{eqnarray*}\label{r7}
        2\bar{a}_{1,2} + 2\bar{a}_{2,1} + 2\bar{a}_{2,2} & = & 0, \\
        2\bar{a}_{1,1} + 2\bar{a}_{2,1} + 2\bar{a}_{2,2} & = & 0, \\
        2\bar{a}_{1,1} + 2\bar{a}_{1,2} + 2\bar{a}_{2,2} & = & 0 \quad \mbox{and} \\
        2\bar{a}_{1,1} + 2\bar{a}_{1,2} + 2\bar{a}_{2,1} & = & 0.
\end{eqnarray*}

We consider the presentation of the homology group of tiling using
the following generators  $x=\bar{a}_{1,2}+\bar{a}_{2,1} + \bar{a}_{2,2}$, $y=\bar{a}_{1,1}+\bar{a}_{2,1} + \bar{a}_{2,2}$, $z=\bar{a}_{1,1}+\bar{a}_{1,2} + \bar{a}_{2,2}$ and
$t=\bar{a}_{1,1} + \bar{a}_{1,2} + \bar{a}_{2,1}+\bar{a}_{2,2}$. The upper
relations in new generators are $$2 x= 2y= 2 z=6 t-2x-2y-2 z=0
.$$ Finally, we find that the homology group of
tiling is
\begin{eqnarray*}
\dst G \langle x, y, z, t
|2c=2b=2z=6t=0 \rangle \cong
\left(\mathbb{Z}_{2}\right)^3\oplus \mathbb{Z}_{6}.
\end{eqnarray*}

\noindent  It follows that the element corresponding to this grid
\begin{eqnarray*}
\Theta & = & (2k+1)\bar{a}_{1,1} + (2k+1)\bar{a}_{1,2}+
(2k+1)\bar{a}_{2,1} + (2k+1)\bar{a}_{2,2} \\ & = &
(2k+1)(\bar{a}_{1,1} + \bar{a}_{1,2} + \bar{a}_{2,1} + \bar{a}_{2,2})
\\ &
= & (2k+1)u
\end{eqnarray*}
is a nontrivial element in the homology group as $6\nmid 2k+1$. Therefore, tiling
does not exist. \qed

\begin{remark}
The same conclusion can be obtained by colouring of torus grid as
in Figure \ref{Teorema_3_3}. Each tile covers $2$ blue cells, $2$
red and $2$  green or  $2$ blue, $2$ yellow and $2$ green or $2$
red, $2$ yellow and $2$ green  or $2$ blue, $2$ red and $2$ green
cells. Given that the number of cells of each color is odd, and
every tile covers even number of cells of each color, we  conclude
that tilling is not possible.
\end{remark}

Now we prove some results on surfaces with boundaries. As we will
see, topology contributes significantly to the homology group of
tiling.

\begin{figure}[!ht]
\centering
    \includegraphics[width=0.9\textwidth]{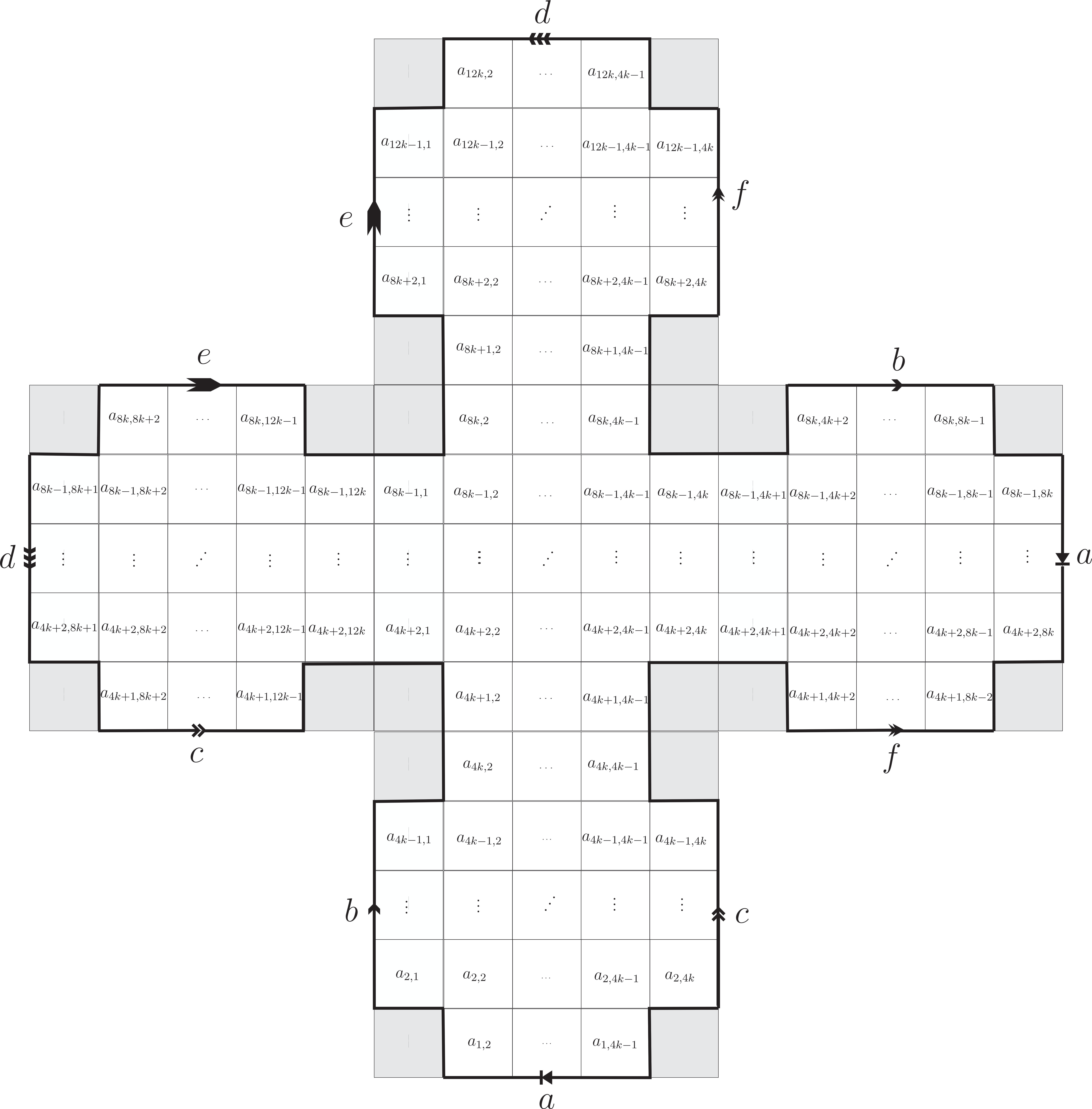}
    \caption{Square grid on a non-orientable surface  of genus $6$ with boundary }\label{I_torus_3}
\end{figure}

\begin{theorem}
A square grid on a non-orientable surface of genus $6$ with
boundary formed by identifying the sides of a dodecagon
consisting of five $4k \times 4k$ squares and removing 20 corner cells around cone point as in  Figure \ref{I_torus_3} cannot be tiled with I-tetrominoes and
Z-tetrominoes.
\end{theorem}

\noindent \dz Let us denote the cells of this square grid as in
Figure \ref{I_torus_3}. Observe that cells $a_{1,1}$, $a_{1,4k}$,
$a_{4k,1}$, $a_{4k,4k}$, $a_{4k+1,8k+1}$, $a_{4k+1,12k}$,
$a_{4k+1,1}$, $a_{4k+1,4k}$, $a_{4k+1,4k+1}$, $a_{4k+1,8k}$,
$a_{8k,8k+1}$, $a_{8k,12k}$, $a_{8k,1}$, $a_{8k,4k}$,
$a_{8k,4k+1}$, $a_{8k,8k}$, $a_{8k+1,1}$, $a_{8k+1,4k}$,
$a_{12k,1}$ and $a_{12k,4k}$ are deleted and that, topologically,
after gluing their union becomes a disk. Thus, we study a gluing
of non-orientable surface  of genus $6$ with one boundary
component.

Using I-tetrominoes it is easy to deduce that in the homology
group of tiling it holds that $\bar{a}_{i, j}=\bar{a}_{i+4, j}$
and $\bar{a}_{i, j}=\bar{a}_{i+4, j}$.

A placement of a Z-tetromino yields one of the following two
relations
$$\bar{a}_{i, j} + \bar{a}_{i+1,j} + \bar{a}_{i+1,j+1} + \bar{a}_{i+2, j+1} =  0 \quad \mbox{and} \quad
\bar{a}_{i+1, j} + \bar{a}_{i+1,j+1} + \bar{a}_{i+2,j+1} +
\bar{a}_{i+2, j+2} =  0.$$ They imply $\bar{a}_{i+2,
j+2}=\bar{a}_{i, j}$.

\begin{figure}[h!]
    \centering
    \includegraphics[width=0.9\textwidth]{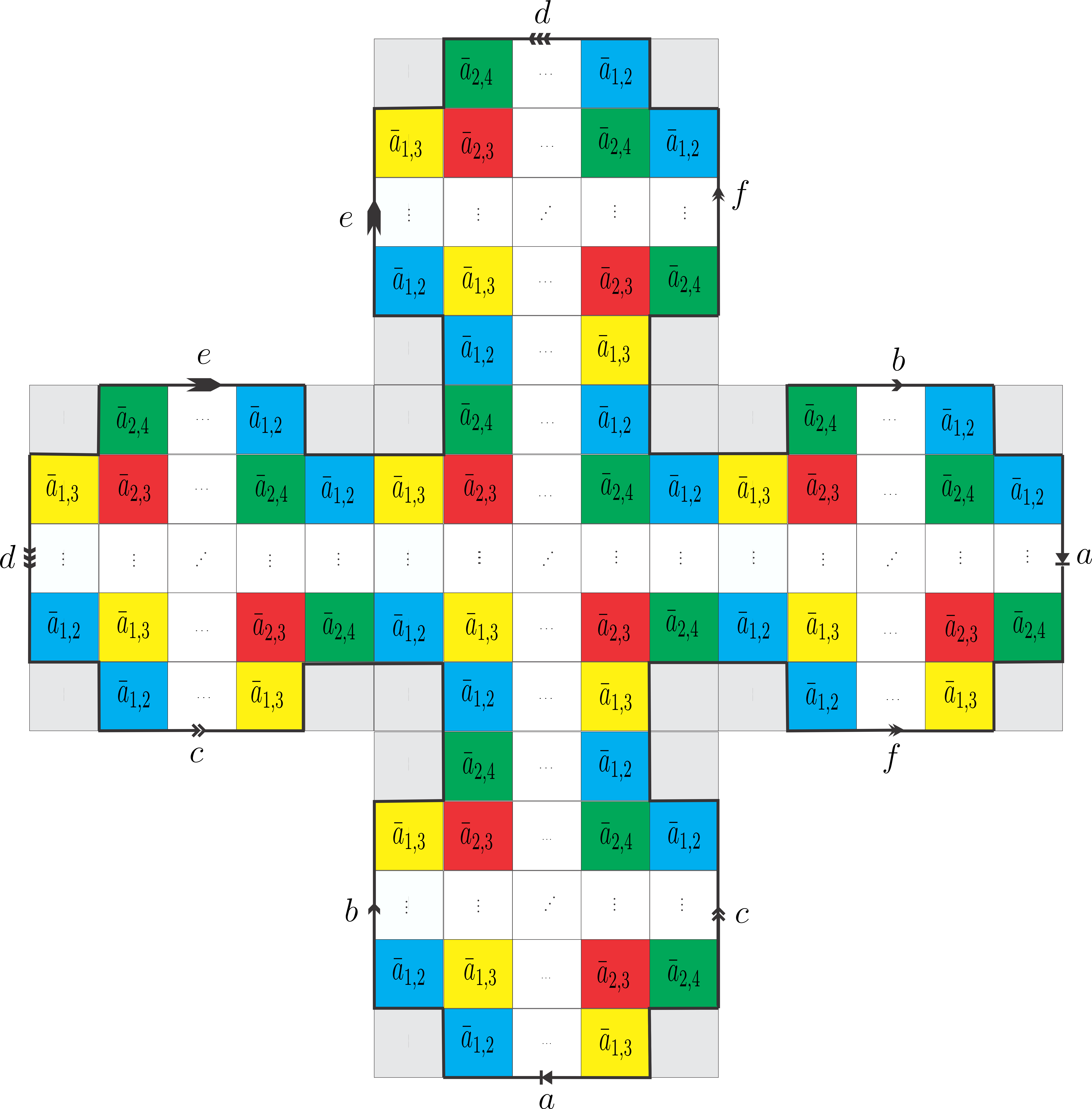}
    \caption{Coloring of the equivalent cells in square grid on a non-orientable surface of genus 6 with boundary}\label{I_torus_3_bojanje}
\end{figure}

Considering placement of I-tetromino across the edge $d$ it is
easy to see that $\bar{a}_{4k+2, 8k+1}=\bar{a}_{12k-3, 2}$,
$\bar{a}_{4k+2, 8k+2}=\bar{a}_{12k-2, 2}$, $\bar{a}_{4k+2,
8k+3}=\bar{a}_{12k-1, 2}$ and $\bar{a}_{4k+2, 8k+4}=\bar{a}_{12k,
2}$. With the relations above we obtain the following equivalences
in the homology group of this tiling depicted in Figure
\ref{I_torus_3_bojanje}.

Thus, the homology group of tiling is a free abelian group with
four generators $\bar{a}_{1,2}, \bar{a}_{1,3}, \bar{a}_{2,3},
\bar{a}_{2,4}$ quotiented by  the following relations
\begin{eqnarray*}
\bar{a}_{1,2}+\bar{a}_{1,3}+\bar{a}_{2,3}+\bar{a}_{2,4}=0,\\
2 \bar{a}_{1,2}+2\bar{a}_{1,3}=0,\\
2 \bar{a}_{1,3}+2\bar{a}_{2,3}=0,\\
2 \bar{a}_{2,3}+2\bar{a}_{2,4}=0,\\
2 \bar{a}_{1,2}+2\bar{a}_{2,4}=0.
\end{eqnarray*}

We eliminate generator $\bar{a}_{2,4}$ from its presentation and
consider generators $ \bar{a}_{1,3}$,
$b=\bar{a}_{1,2}+\bar{a}_{1,3}$ and
$c=\bar{a}_{1,3}+\bar{a}_{2,3}$. We obtain that our group of
homology is isomorphic to $\mathbb{Z}\oplus\mathbb{Z}_2^2$.

Our square grid contains
 $20k^2$ cells $\bar{a}_{1,2}$, $5(4k^2-1)$ cells $\bar{a}_{1,3}$ and $\bar{a}_{2,4}$, as well as $10(2k^2-1)$ cells $\bar{a}_{2,3}$. The element corresponding to this grid
\begin{eqnarray*}
\Theta & = & 20k^2\bar{a}_{1,2} + 5(4k^2-1)\bar{a}_{1,3} + 10(2k^2-1)\bar{a}_{2,3}+5(4k^2-1)\bar{a}_{2,4}\\
& = & 20k^2(\bar{a}_{1,2}+\bar{a}_{1,3}+\bar{a}_{2,3}+\bar{a}_{2,4}) -5\bar{a}_{1,3} - 10\bar{a}_{2,3} -5\bar{a}_{2,4}\\
& = & 5 \bar{a}_{1,2}-5\bar{a}_{1,3} =b-10\bar{a}_{1,3}
\end{eqnarray*}
is a non-trivial element of the  homology group and desired tiling
is not possible. \qed

\begin{figure}[h!]
    \centering
    \includegraphics[width=0.85\textwidth=0.95]{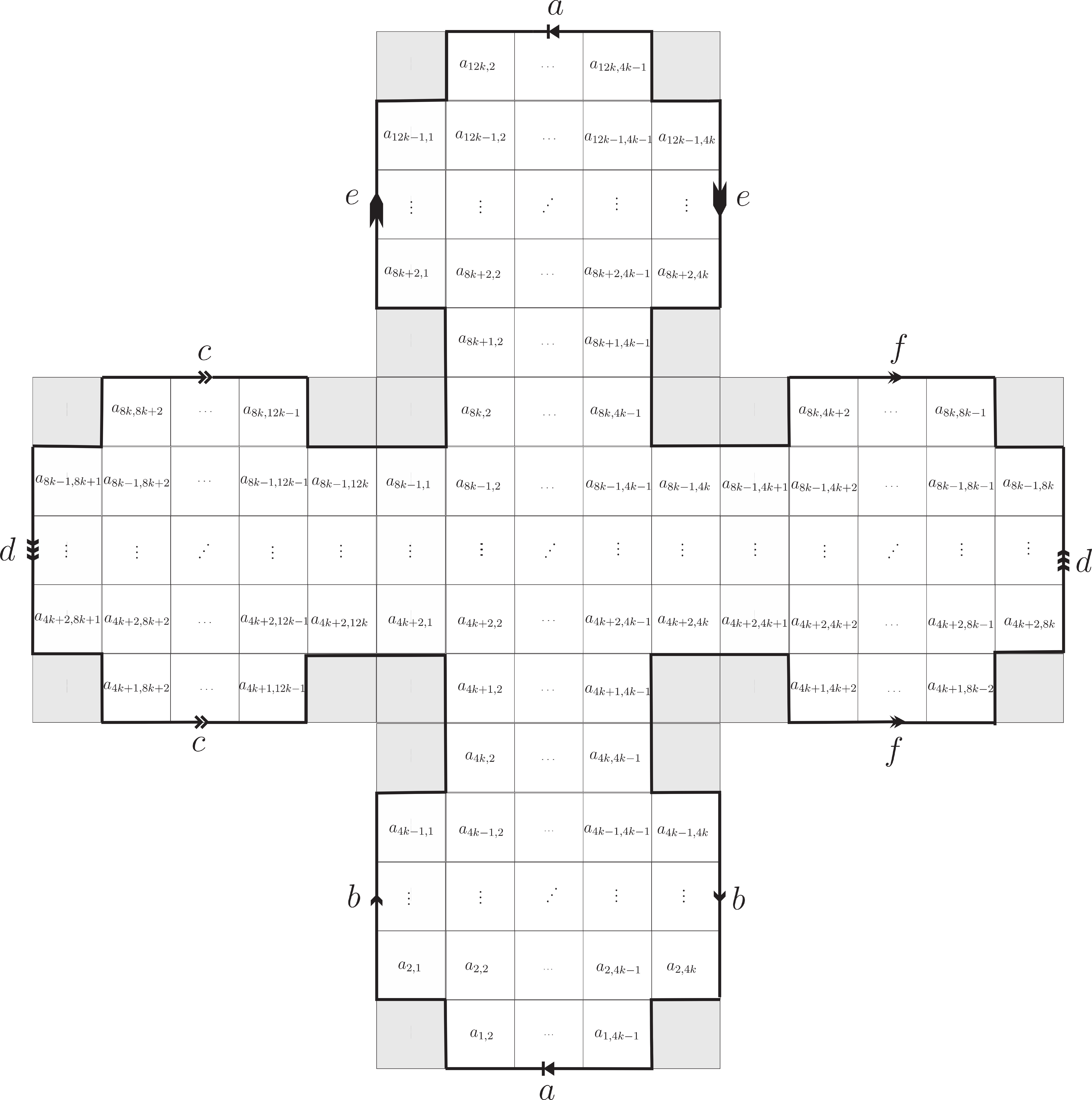}
    \caption{Square grid on a non-orientable surface of genus $4$ with three boundary components}\label{L_torus_3}
\end{figure}

\begin{theorem}
A grid on a non-orientable surface of genus $4$ with boundary is
formed by identifying the sides of a dodecagon consisting of five
$4k \times 4k$  squares and with removed $20$ cells around cone points as in Figure
\ref{L_torus_3} cannot be tiled with L-tetrominoes.
\end{theorem}

\noindent \dz Model in Figure \ref{L_torus_3} after gluing along
marked sides and deletion of $20$ corner cells gives a
non-orientable surface of genus $4$ with three boundary
components. Denote the cells in the grid as in the previous
example.

Placing  L-tetromino in the given model in vertical position
before taking identification into account will give one of the two
relations
\begin{eqnarray}\label{vertikalL1}
\bar{a}_{i, j} + \bar{a}_{i+1,j} + \bar{a}_{i+2,j} + \bar{a}_{i+2,
j-1}& = & 0 \quad \mbox{and} \\ \label{vertikalL2} \bar{a}_{i, j}
+ \bar{a}_{i+1,j} + \bar{a}_{i+2,j} + \bar{a}_{i+2, j+1}& = & 0
\end{eqnarray}
in the homology group of tiling. From  \eqref{vertikalL1} and
\eqref{vertikalL2} we obtain that in the group of homology of this
tiling the cells $\bar{a}_{i,j-1}=\bar{a}_{i,j+1}$ are equivalent.
Analogously, it holds  that $\bar{a}_{i-1,j}=\bar{a}_{i+1,j}$ are
equivalent in the homology group of this tiling.

 We summarize all upper equivalences of cells   in
\begin{eqnarray*}
\bar{a}_{i, j}=\begin{cases}
\bar{a}_{1, 1}, \quad \text{if} \quad i \equiv 1 \pmod{2}, \, j \equiv 1 \pmod{2},\\
\bar{a}_{1, 2}, \quad \text{if} \quad i \equiv 1 \pmod{2}, \, j \equiv 0 \pmod{2},\\
\bar{a}_{2, 1}, \quad \text{if} \quad i \equiv 0 \pmod{2}, \, j \equiv 1 \pmod{2},\\
\bar{a}_{2, 2}, \quad \text{if}  \quad i \equiv 0 \pmod{2}, \, j \equiv 0 \pmod{2}.\\
\end{cases}
\end{eqnarray*}

\begin{figure}[h!]
    \centering
    \includegraphics[width=0.65\textwidth=0.95]{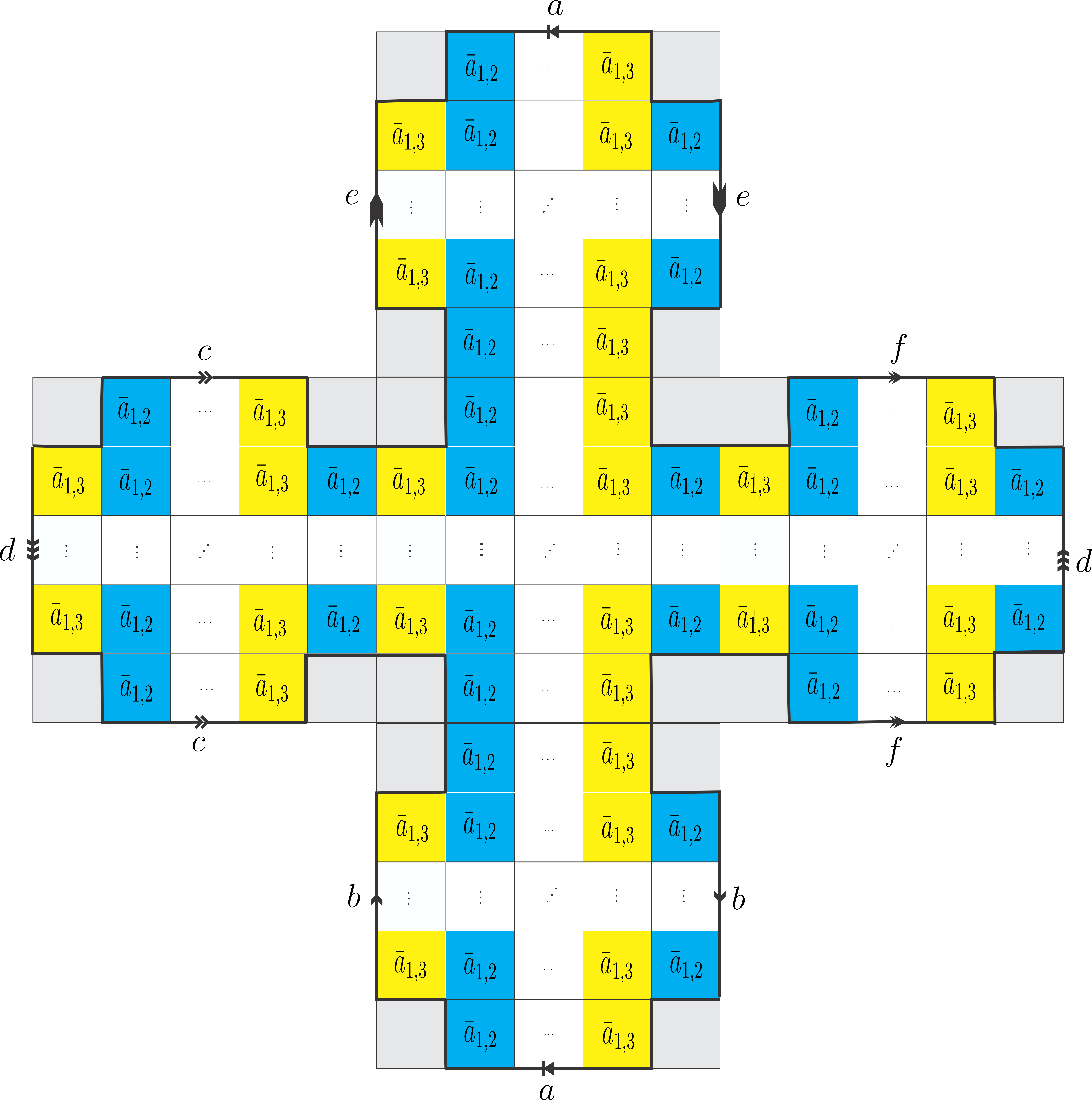}
    \caption{Equivalent cells on square grid on a non-orientable surface of genus $4$  with boundary}\label{L_torus_3_bojanje}
\end{figure}

Consider a placement of L-tetromino along edge  denoted by $e$ in
Figure \ref{L_torus_3} and corresponding equations in the homology
group of tiling
\begin{eqnarray*}\label{e}
\bar{a}_{1, 1} + \bar{a}_{1,2} + \bar{a}_{2,1} + \bar{a}_{1, 2}& = & 0 \quad \mbox{and} \\
\bar{a}_{2, 1} + \bar{a}_{1,2} + \bar{a}_{2,1} + \bar{a}_{1, 2}& =
& 0.
\end{eqnarray*} From them we deduce that $\bar{a}_{1, 1}=\bar{a}_{2, 1}$. In a similar way we obtain that  $\bar{a}_{1, 2}=\bar{a}_{2, 2}$.  These equivalences are illustrated in Figure \ref{L_torus_3_bojanje}.

Placement of $L$ tetromino on the grid with equivalent cells,
including placements across glued sides, we obtain one of the two
relations
\begin{eqnarray*}
3\bar{a}_{1,1}+\bar{a}_{1,2} & = & 0 \quad \mbox{and}\\
3\bar{a}_{1,2}+\bar{a}_{1,1} & = & 0.
\end{eqnarray*}

Now we conclude that $8 \bar{a}_{1, 1}=0$. Therefore, the homology
group is isomorphic to the group
\begin{eqnarray*}
\dst G \langle \bar{a}_{1,1}| 8 \bar{a}_{1, 1}=0 \rangle\cong
\mathbb{Z}_8.
\end{eqnarray*}
Our square grid contains $10(4k^2-1)$ cells $a_{1,1}$ and
$a_{1,2}$, so the element assigned to this grid
\begin{eqnarray*}
\Theta & = & 10(4k^2-1)\bar{a}_{1,1} + 10(4k^2-1)\bar{a}_{1,2} = 4
\bar{a}_{1, 1}
\end{eqnarray*}
is a non-trivial element in the homology group of tiling and it is
not possible to tile the given grid using L-tetrominoes.  \qed

\begin{figure}[h!]
\centering
    \includegraphics[width=0.65\textwidth]{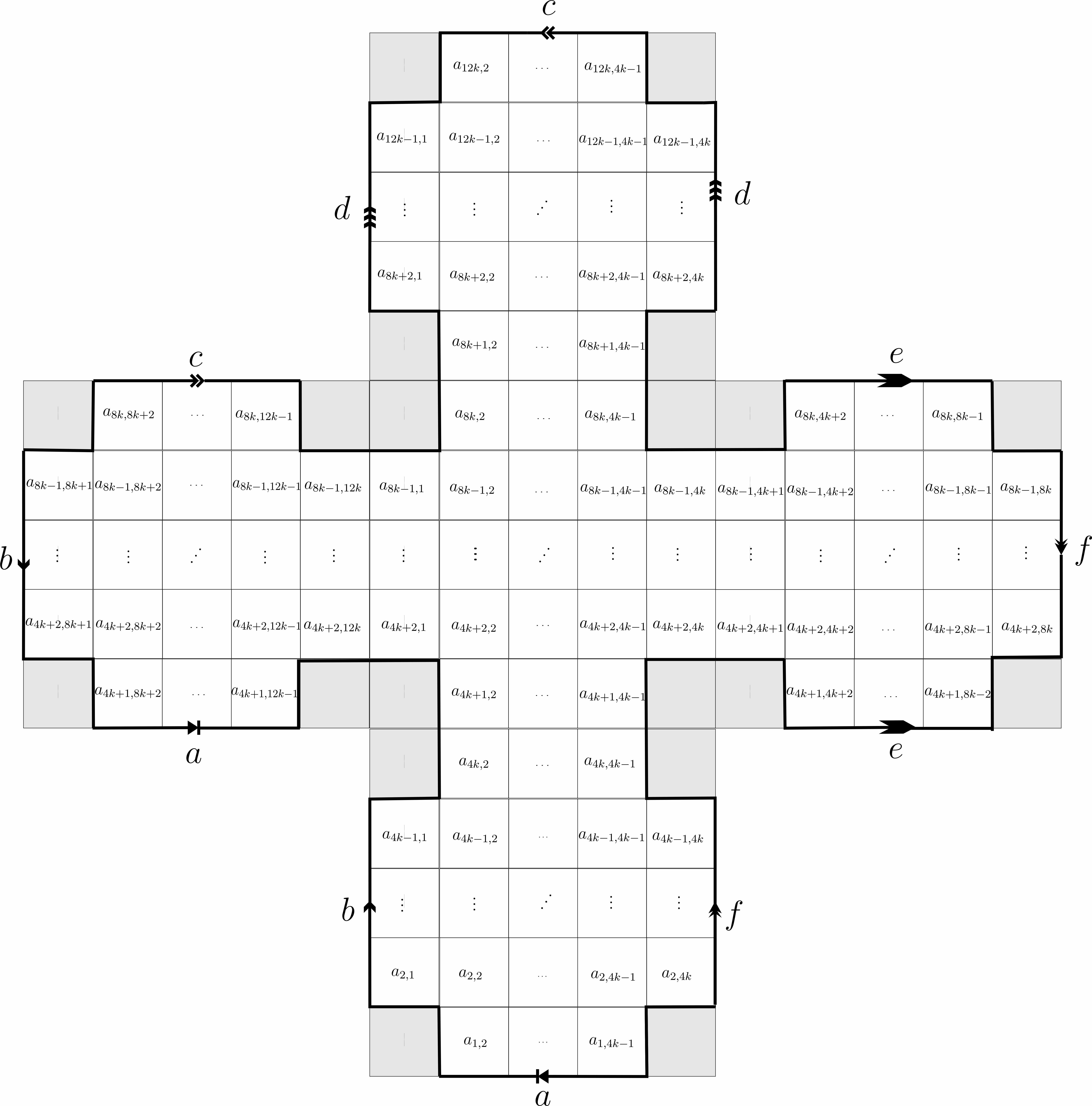}
    \caption{Square grid on an orientable genus $3$ surface with boundary}\label{T_torus_3}
\end{figure}

\begin{theorem}
A square grid on an orientable surface of genus $3$ with boundary
formed by identifying the sides of a dodecagon  consisting of five
$4k \times 4k$ squares and removing 20 cells meeting in the cone point as in Figure
\ref{T_torus_3} cannot be tiled by  $T$-tetrominoes.
\end{theorem}

\noindent \dz

\begin{figure}[H]
    \centering
    \includegraphics[width=0.75\textwidth]{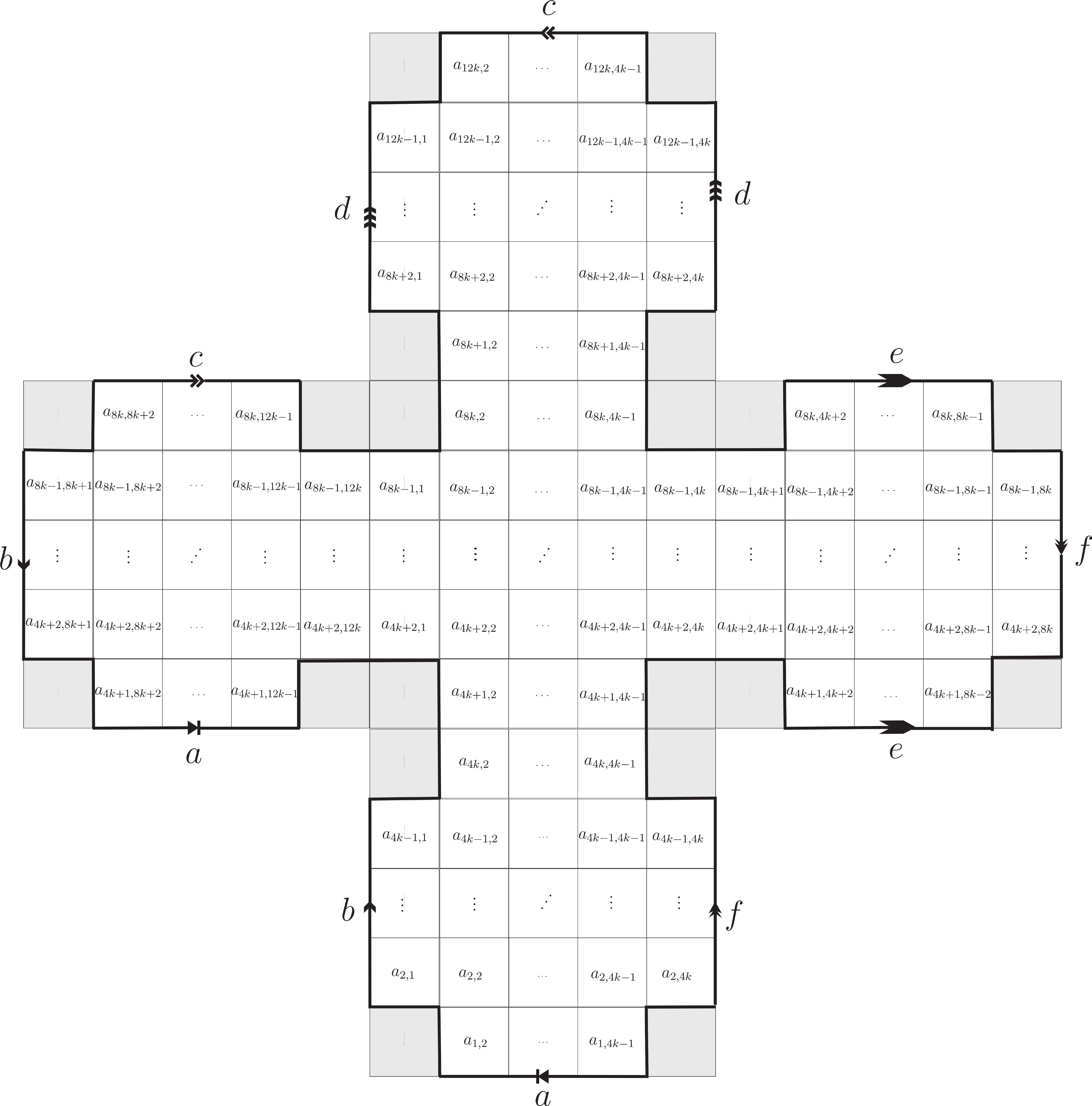}
    \caption{Equivalent cells on square grid on an orientable genus $3$ surface with boundary}\label{T_torus_bojanje}
\end{figure}

It is straightforward to check that model in Figure
\ref{T_torus_3} after gluing along marked sides and deletion of 20
corner cells gives a genus $3$ surface with one boundary
component. Denote the cells in the grid as in the previous
theorem. The following equality is easily obtained

 \begin{eqnarray*}
\bar{a}_{i, j}=\begin{cases}
\bar{a}_{1, 1}, \quad \text{if} \quad i -j\equiv 0 \pmod{2}, \\
\bar{a}_{1, 2}, \quad \text{if} \quad i-j \equiv 1 \pmod{2},
\end{cases}
\end{eqnarray*}
as it is illustrated in Figure \ref{T_torus_bojanje}.

\noindent If we put $T$ tetrominoes on the grid with equivalent
cells, even placing it across a glued sides, we obtain one of the
two relations
\begin{eqnarray*}
3\bar{a}_{1,3}+\bar{a}_{1,2} & = & 0 \quad \mbox{and} \\
3\bar{a}_{1,2}+\bar{a}_{1,3} & = & 0.
\end{eqnarray*}
Therefore, we get that the homology group of this is isomorphic to
the group
\begin{eqnarray*}
\dst G \langle \bar{a}_{1,2} | 8\bar{a}_{1,2}=0 \rangle \cong
\mathbb{Z}_{8}.
\end{eqnarray*}

Our square grid contains $10(4k^2-1)$ cells $\bar{a}_{1,2}$ and
$10(4k^2-1)$ cells $\bar{a}_{1,3}$, so the element assigned to
this grid is
\begin{eqnarray*}
\Theta & = & 10(4k^2-1)\bar{a}_{1,2} + 10(4k^2-1)\bar{a}_{1,3} \\
& = &   40k^2\bar{a}_{1,2} - 10\bar{a}_{1,2}  -120 k^2
\bar{a}_{1,2}+30 \bar{a}_{1,2}  = 4 \bar{a}_{1,2}.
\end{eqnarray*}
$\Theta$ is  a non trivial element in the homology group of
tiling, and therefore it is not possible to til the given square
grid using T-tetrominoes. \qed

\begin{figure}[!ht]
\centering
    \includegraphics[width=\textwidth]{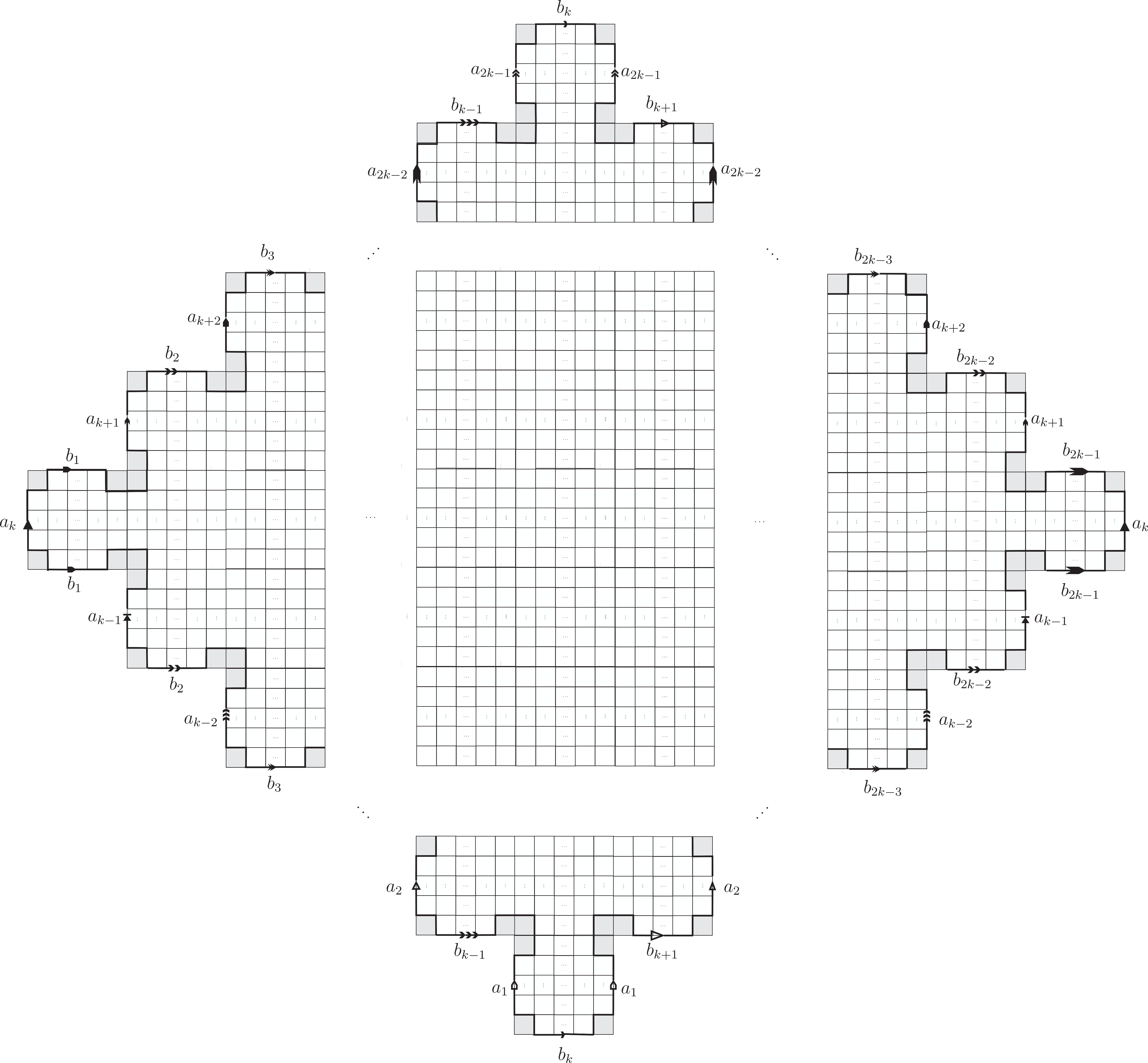}
    \caption{Square grid on an orientable genus $2k-1$ surface with boundary}\label{grand}
\end{figure}

\begin{theorem}\label{teo3.7}
A square grid on an orientable surface of genus $2k-1$ with
boundary formed by identifying the sides of a $(8k-4)$-gon
consisting of $2 k^2-2 k+1$ squares of side $(4k-3)d$ where $d$ is
a positive integer,  without corner cells as in Figure \ref{grand}
can not be tiled with $1\times (4k-3)$ polyomino.
\end{theorem}

\noindent \dz  From Figure \ref{grand} it is clear that the
surface is orientable. Label the cells in the grid in standard
way. Denote the cell by $a_{i, j}$ in standard way assuming that
the bottom left corner cell is $a_{1, 1}$. As with other
$I$-minoes it is straightforward to get
$$\bar{a}_{i,j}=\bar{a}_{i,j+4k-3}=\bar{a}_{i+4k-3,j}.$$

Using this equivalences we find that there are $(4k-3)^2$ types of
the cells $\bar{a}_{i, j}$, $1\leq i, j\leq 4k-3$ in the homology
group of tiling. We see that there are $8k-6$ relations
\begin{eqnarray*}
\sum_{j=1}^{4k-3} \bar{a}_{i, j} =0 \quad \mbox{for} \, i=1, \dots, 4 k-3 \quad \mbox{and}\\
\sum_{i=1}^{4k-3} \bar{a}_{i, j} =0 \quad \mbox{for} \, j=1,
\dots, 4 k-3
\end{eqnarray*} assigned to a placement of $1\times (4k-3)$ polyomino on the board (including placements across gluing sides). Therefore, our homology group of tiling is isomorphic to
$$G \langle a_{i,j}| 1\leq i, j\leq 4k-4 \rangle\cong \mathbb{Z}^{16(k-1)^2}.$$

Element $\Theta$ assigned to the grid is $$\Theta=-(2k-1)
\sum_{i=2}^{4k-4} \sum_{j=2}^{4k-4} \bar{a}_{i, j}.$$ This is a
non-trivial element in the homology group of tiling and the claim
is therefore proved. \qed

\bigskip

\begin{center}\textmd{\textbf{Acknowledgements} }
\end{center}

\medskip The authors are grateful to Djordje \v{Z}ikeli\'{c} and Igor Spasojevi\'{c} for valuable comments and discussions. The second author was supported by the Ministry for Education, Science and Technological Development of the Republic of Serbia through the Mathematical Institute SANU.

\end{document}